\documentstyle[12pt]{article}

\advance\hoffset by -7mm  
\makeatletter
\@addtoreset{equation}{section}
\makeatother

\input amssym.def  
\input amssym.tex

\renewcommand{\title}[1]{\null\vspace{25mm}

\noindent{\Large{\bf #1}}\vspace{10mm}

\noindent {\large By }}
\newcommand{\authors}[1]{\noindent{\large #1}\vspace{3mm}

}
\newcommand{\address}[1]{\noindent #1\vspace{5mm}

}
\renewcommand{\abstract}[1]{\vspace{19mm}

\noindent{\small{\em Abstract.} #1}\vspace{2mm}

}

\begin{document}

\title{Rigidity for periodic magnetic fields}
\authors{M L Bialy}
\address{School of Math Sciences, Tel-Aviv University, Tel-Aviv 69978,
  Israel\footnote{The author was partially supported by EPSRC and
    Arc-en-Ciel/98.}\\
Email:  bialy@math.tau.ac.il}

\abstract{
We study the motion of a charge on a conformally flat 
Riemannian torus in the presence of magnetic
field.  We prove that for any non-zero magnetic field 
there always exist orbits of this motion which have conjugate
points.  We conjecture that the restriction of conformal flatness 
of the metric is
not essential for this result.  This would provide a ``twisted''
version of the recent generalisation of Hopf's rigidity result
obtained by Burago and Ivanov.
}

\section{Introduction}

It was proved by E. Hopf \cite{[H]}, that for Riemannian 2-torus
with non-zero curvature there always exist geodesics with conjugate
points.  In \cite{[B-I]} this result was finaly generalised to higher
dimensions.  We refer the reader to \cite{[K]}, \cite{[C-K]},
\cite{[C-F]} and \cite{[He]} for a non-complete list of previous
important contributions.  It was discovered in \cite{[B]},
\cite{[B-P1]} ( see also \cite{[B-M]}) that the nature of E. Hopf
rigidity is not entirely Riemannian and can be established for other
variational problems.  In each case this requires some new integral
geometric tool adapted for the system under consideration.

In this paper we give a proof of such a rigidity result for the motion
of a charged particle on a torus in the presence of the magnetic field, 
provided the
Riemannian metric is conformally flat.  This result gives the twisted version
of a theorem by A. Knauf and C .Croke-A. Fathi (\cite{[C-F]} and \cite{[K]}).  
We formulate the conjecture that the
restriction of conformal flatness is not essential for the result.
This, being true, would provide a twisted version of Burago-Ivanov
theorem.

In the situation we consider, it would be reasonable to think that
there always exist a periodic orbit with conjugate points.  This fact
is, however, much more difficult to approach.  It is not proved to the
best of our knowledge even in the usual untwisted case.  
Moreover, it seems that the
fact of existence of periodic orbit is still not completely answered
for arbitrary magnetic field on the torus.  We refer the reader to a survey 
paper by V. Ginzburg
\cite{[Gi]} for known results and techniques.

\section{Main Results}

Magnetic field on a torus ${\Bbb T}^n$ is a closed 2-form which we will
denote by $\beta$.  Consider the cotangent bundle $T^* {\Bbb T}^n$ and
define the symplectic structure twisted by $\beta$.
$$
\omega=\omega_0 + \pi^* \beta
$$
where $\omega_0$ is a standard structure and $\pi$ is the canonical
projection.

For a Riemannian metric $g$ on ${\Bbb T}^n$ consider the Hamiltonian
flow $g^t$ of the function $H$
$$
H={1 \over 2} <p,p>_g
$$
computed with the help of the symplectic structure $\omega$.  Let me
remind that no conjugate points condition for the orbit $\vartheta$ of
$g^t$ means that for any two points $x, y \in \vartheta$
$$
g^t_* V(x)\pitchfork V(y)
$$
where $t$ is a time difference between $y$ and $x$ and $V(x)$, $V(y)$
are vertical subspaces at $x,y$.

\newtheorem{theorem}{Theorem} 
\begin{theorem}
\label{T21}
Suppose that the Riemannian metric $g$ is conformally flat.  Then
there always exist orbits of $g^t$ on the level of $\{ H= 1/2 \}$ with
conjugate points unless the 2-form $\beta$ vanishes identically.
\end{theorem}

\newtheorem{remark}{Remark}
\begin{remark}
\label{R22}
In other words, no conjugate points condition implies $\beta \equiv 0$ and then
it follows from \cite{[C-F]} and \cite{[K]} that the metric is flat too.
Actually, we will see it once more by our computations later.
\end{remark}

The proof which is suggested in this paper follows the 
original scheme by E. Hopf.  However, it
might need serious modifications if one tries to generalize the result for
arbitrary Riemannian metric.

The first ingredient of the proof is to construct measurable field of
Lagrangian subspaces $l(x) \subseteq T_x (T^* {\Bbb T}^n)$ for any $x
\in \{ H= 1/2 \}$.  This field can be constructed by the
following limit procedure:
$$
l(x) = lim_{t \rightarrow + \infty} g^t_* V \left( g^{-t} (x) \right)
$$

It was first used by E. Hopf \cite{[H]} and L. Green \cite{[G]} for
Riemannian case.  We refer the reader to recent paper \cite{[C-I]} for
the proof in a general optical case.

It follows from the very construction of $l$ that the field $l$ is
invariant under $g^t$ and is transversal to the vertical field $V$
everywhere.

With the construction of $l$ Theorem \ref{T21} is a corollary of
the following 

\begin{theorem}
\label{T23}
Let $l$ be a measurable field of Lagragian subspaces invariant under
the flow $g^t$.  If $l$ is transversal to $V$ everywhere then $\beta$
vanishes identically, and the metric $g$ is flat.
\end{theorem}

This last theorem has the following dynamical interpretation.

\begin{theorem}
\label{T24}
Suppose that the energy shell $\{ H= 1/2 \}$ is smoothly foliated by
Lagrangian tori homologous to the zero section of $T^* {\Bbb
  T}^n$.  Then the 2-form $\beta$ vanishes identically and the metric
is flat.
\end{theorem}

\begin{remark}
\label{R25}
In the paper \cite{[K]}, such a situation is called total
integrability and simple examples of totally integrable magnetic
geodesic flows are given.  Let me remark that the Lagrangian torii in
all these examples cannot be homologous to the zero section as it
follows from Theorem \ref{T24}.  It is an interesting, completely open problem
to characterise totally integrable magnetic geodesic flows.
\end{remark}

\noindent{\large \bf Acknowledgements}

\noindent
I was introduced to the subject of magnetic fields by Viktor Ginzburg
(UC Santa Cruz).  He suggested to me the question on the rigidity for
magnetic fields.  I am deeply grateful to him for many very useful discussions.

The results of this paper were presented at the Symplectic Geometry meeting
in Warwick 1998 and on the Arthur Besse Geometry Seminar.  I am
thankful to D. Salamon, F. Laudenbach and P. Gauduchon for inviting me
to speak there.

I would also like to thank the EPSRC and Arc-en-Ciel for their
support.

\section{Proofs}

Let me explain first how Theorem \ref{T24} follows from Theorem \ref{T23}.
Since the torii are Lagrangian and homologous to the zero section
then the 2-form $\beta$ must be exact.  Denote by $\alpha$ the
primitive 1-form.  Then one can easily see that the flow $g^t$ is
equivalent to the Hamiltonian flow $\tilde g^t$ of the function $\tilde
H = 1/2 < p - \alpha, p - \alpha>_g$ with respect to the standard
structure $\omega_0$.  This equivalence is given by the diffeomorphism
$(q,p) \rightarrow (q,p + \alpha)$ which is fiber preserving.  Note
that the function $\tilde H$ is strictly convex with respect to $p$.
It follows from generalised Birkhoff theorem (see \cite{[B-P2]} for
its most general form and for the survey and discussions) that all the
Lagrangian torii are the sections of the cotangent bundle.  But then the
distribution of their tangent spaces meets the conditions of Theorem
\ref{T23}.

{\bf Proof of Theorem \ref{T23}}

We shall work in standard coordinates $(q,p)$ on $T^* {\Bbb
  T}^n$, such that the Riemannian metric $g$ is given by
$$
ds^2 = {1 \over 2 \lambda} \left( dq^2_1 + \cdots + dq^2_n \right).
$$
Then
\begin{equation}
H(p,q) = {\lambda \over 2} \left( p^2_1 + \cdots + p^2_n \right)
\label{31}
\end{equation}

Write $\beta = d \alpha + \gamma$, where $\alpha$ is a 1-form and
$\gamma$ is a 2-form, having constant coefficients in the coordinates
$(q_1 \cdots q_n)$,
$\gamma = \Sigma_{i<j} \; \gamma_{ij} \; dq_i \; \wedge
\; dq_j$.  By the change of coordinates
$$
(q,p) \rightarrow (q,p + \alpha)
$$
we have an equivalent Hamiltonian flow $\tilde g^t$ of
$$
\tilde H (p,q) = H (p- \alpha, q) = {1 \over 2} \lambda \left( (p_1 -
  \alpha_1)^2 + \cdots +(p_n - \alpha_n)^2 \right)
$$
with respect to the symplectic form
$$
\tilde \omega = \omega_0 + \gamma = \Sigma^n_{i=1} dp_i \wedge dq_i +
\Sigma_{i<\gamma} \gamma_{ij} \; dq_i \; \wedge \; dq_j.
$$ 

Denote by $\tilde l$ the invariant distribution of Lagrangian
subspaces with respect to $\tilde \omega$.  At any point $x=(q,p)
\in \{\tilde H=1/2 \}$ it is given by a matrix $A(x)$:
$$
dp = A(p,q)dq
$$

The condition to be Lagrangian is equivalent to
$$
A^T - A = \Gamma
$$
where $\Gamma$ is the skew symmetric matrix of $\gamma$.  Then $A$ is
a measurable matrix function and satisfies the following Ricatti
equation along the flow $\tilde g^t$.
$$
\dot A + (A + \Gamma) \tilde H_{pp} A + (A + \Gamma) \tilde H_{pq} +
\tilde H_{qp} A + \tilde H_{qq} = 0
$$
where $\tilde H_{pp}$, $\tilde H_{pq}$, $\tilde H_{qq}$ are the
matrices of second derivatives of $\tilde H$.  Then it can be written in
the form
\begin{eqnarray*}
\dot A &+& \left(( A + \Gamma) \tilde H_{pp}^{1/2} + \tilde H_{qp} \tilde
  H_{pp}^{- 1/2}\right) \left( \tilde H_{pp}^{1/2} A + \tilde
  H_{pp}^{- 1/2} \tilde H_{pq} \right)\\
&+& \left( \tilde H_{qq} - \tilde H_{qp} \tilde H_{pp}^{-1} \tilde
  H_{pq} \right) = 0
\end{eqnarray*}

Introduce the function $a (p,q) = tr A (p,q)$.  Then it follows
$$
\dot a + tr \left( \tilde H_{qq} - \tilde H_{qp} \tilde H_{pp}^{-1} \tilde
  H_{pq}\right) \leq 0
$$
with equality possible only when $A= - \tilde H_{pp}^{-1} \tilde
  H_{pq}$.  Integrate now this inequality with respect to the
  invariant measure $\tilde \mu$ on the energy level {$\tilde H =
    1/2$}.  We get the following
\begin{equation}
\sigma (\tilde H) = \int tr \left( \tilde H_{qq} - \tilde H_{qp}
  \tilde H_{pp}^{-1} \tilde H_{pq} \right) d \tilde \mu \leq 0
\label{32}
\end{equation}
with equality possible only for $A=-\tilde H_{pp}^{-1} \tilde H_{pq}$
everywhere on the level $\{ \tilde H = 1/2 \}$.  On the other hand, we
shall compute $\sigma (\tilde H)$ by the following:

\newtheorem{lemma}{Lemma} 
\begin{lemma}
\label{L33}
For any Hamiltonian function $H(p,q)$ convex and symmetric in $p$
$H(p,q) = H(-p,q)$ it follows that $\sigma (\tilde H) = \sigma (H)$
where $\tilde H (p,q) = H (p - \alpha (q),q)$ for any 1-form
$\alpha$.  Here $\sigma (H)$ is computed by the integral 
(3.2) with respect to the hamiltonian $H$ and the invariant measure 
$\mu$ of $g^t$.
\end{lemma}

\begin{lemma}
\label{L34}
If $H$ is given by (3.1) then $\sigma (H) \geq 0$ and equality
is achieved only when $\lambda = const$.
\end{lemma}

\begin{remark}
\label{R35}
The last Lemma is actually the inequality proved in \cite{[K]}, but
here it appears in the Hamiltonian version.  
\end{remark}

Let me complete the proof of the theorem, postponing ones of the
Lemmas.

Comparing (\ref{32}) with Lemmas \ref{L33} and \ref{L34}, one
concludes that the metric is Euclidean and 
$$
A=-\tilde H_{pp}^{-1} \tilde H_{pq}, \quad \mbox{where} \quad \tilde H
= {1 \over 2} \Sigma^n_{i=1} \left( p_i - \alpha_i \right)^2
$$
Then
\begin{equation}
A = \left( A_{ij} \right) = - \left( \tilde H_{p_i q_j} \right) = \left(
  {\partial \alpha _i \over \partial q_j}\right)
\label{36}
\end{equation}

Compute and plug the second derivatives in the Ricatti equation. We have
$$
(t r \dot A) = -tr \left( \tilde H_{qq} - \tilde H_{qp} \tilde
  H_{pp}^{-1} \tilde H_{pq} \right) = \Sigma^n_{ij=1} \left( p_j -
  \alpha_j \right) {\partial^2 \alpha_j \over \partial q^2_i}
$$

Differentiating $trA$ explicitly along the flow $\tilde g^t$, one gets
another expression
$$
(t r \dot A) = \Sigma^n_{j=1} {\partial (trA) \over \partial q_j} \tilde
H_{p_j} = \Sigma^n_{ij=1} (p_j - \alpha_j) {\partial^2 \alpha_i \over
  \partial q_i \; \partial q_j}
$$

Comparing these two expression we obtain
$$
{\partial (trA) \over \partial q_j} = \Delta \alpha_j \quad \mbox{for}
\quad j=1, \dots n.
$$

Then the compatibility condition gives
$$
\Delta \left( {\partial \alpha_k \over \partial q_j} - {\partial \alpha
   _j \over \partial q_k} \right) = 0 \quad \mbox{for any} \quad j,
k=1, \dots n. 
$$

This implies
$$
{\partial \alpha_k \over \partial q_j} - {\partial \alpha_j \over
  \partial q_k} \equiv const = 0 \quad \mbox{and thus} \quad d \alpha
=0.
$$

Also by (\ref{36}) we have $A^T - A = \Gamma = 0$.  So $\beta = d
\alpha + \gamma = 0$.  This completes the proof of the theorem. $\Box$

{\bf Proof of Lemma {\ref{L33}}}

We have to show that $\sigma$ is the same for $H$ and $\tilde H = H (p
- \alpha(q),q)$.  The direct computation gives the following
expressions for the matrices of second derivatives.

\begin{eqnarray*}
\tilde H_{pp} &=& H_{pp}\\
\tilde H_{pq} &=& -H_{pp} D \alpha + H_{pq}, \quad \tilde H_{qp} = -
(D \alpha)^T H_{pp} + H_{qp}\\
\tilde H_{qq} &=& (D \alpha)^T H_{pp} (D \alpha) - H_{qp} D \alpha -
(D \alpha)^T H_{pq} + H_{qq} + M\\
\mbox{where}\quad M &=& (M_{ij}), \quad M_{ij} = \Sigma^n_{k=1} H_{p_k}
{\partial^2 \alpha_k \over \partial q_i \; \partial q_j}
\end{eqnarray*}

Thus we have after performing change of variables 
$(q,p) \rightarrow (q,p - \alpha)$
$$
\sigma (\tilde H) = \sigma (H) + \int_{ \left\{ H = {1 \over 2} \right\} }
\Sigma^n_{k,i=1} H_{p_k} {\partial^2 \alpha_k \over \partial q^2_i} d \mu
$$

By the symmetry assumption this additional integral vanishes because
the integrand is odd function and the level is even.  This completes
the proof of the lemma. $\Box$ 

{\bf Proof of Lemma {\ref{L34}}}

The second derivative for the Hamiltonian $H={\lambda/2} \left( p^2_1
  + \cdots + p^2_n \right)$ are given by the formulas
$$
H_{p_ip_j} = \lambda \delta_{ij}, \quad H_{q_iq_j} = {1 \over 2}
(\lambda)_{q_iq_j} \left( p^2_1 + \cdots + p^2_n \right), \quad H_{p_iq_j}
= \lambda_{q_j} p_i
$$

Thus for the points of the energy level $\left\{ H = {1 \over 2} \right\}$
$$
tr \left( H_{qq} - H_{qp} H^{-1}_{pp} H_{pq} \right) = {1 \over 2 \lambda} \Delta \lambda - {1 \over 2 \lambda^2} (grad \lambda)^2
$$

The invariant measure $\mu$ on $\left\{ H = {1 \over 2} \right\}$ is
computed from the condition
$$
d \mu \wedge d H = \omega^n_0
$$
and then can be easily computed and equals
$$
d \mu = \left( {1 \over \sqrt \lambda} \right)^n d \sigma dq
$$
where $d \sigma$ is a standard measure on unit sphere in $p$ space.
Then we have 
$$
\sigma (H) = Vol (S^{n-1}) \int \left( {1 \over 2 \lambda} \Delta \lambda -
  {1 \over 2 \lambda^2} (grad \lambda)^2 \right) \left( 1 \over \sqrt
  \lambda \right)^n dq
$$

Integrating by parts we obtain
$$
\sigma (H) = Vol (S^{n-1}) \left( {n-2 \over 4} \right) \int
\lambda^{-2 -{n \over 2}} (grad \lambda)^2 dq
$$

From this expression the assertion follows.  $\Box$ 

\section{Some Open Problems}

Let me discuss here some natural open problems.

In the proof of Theorem \ref{T24} we noticed first that the form
$\beta$ must be exact and then showed that it vanishes.  It may happen
that the condition of existence of invariant distribution of
Lagrangian non-vertical subspaces also implies that the form $\beta$
is exact.  We were not able, however, to find a direct proof of this
fact.

It is an important problem to study those magnetic fields having
integrable twisted geodesic flow.  The simplest possible topology of
the phase portrait is one for the totally integrable flows.  It is not
clear how to characterise those flows even for dimension 2.  One
natural class of examples is for $S^1$ symmetric magnetic field $\beta
= \beta (q_1) dq_1 \wedge dq_2$.  We don't know if there are other
examples.  Note that the Lagrangian torii of the foliation can not be
homologous to the zero section ( theorem 3) .

One would like to generalise our results to the case of any
Riemannian metric on ${\Bbb T}^n$.  However, it may happen that
twisted version of Burago-Ivanov proof should be found.
It should be mentioned that our approach works each time when
it is related to a good choice of coordinates in the phase space.
It is not clear to me however, if the coordinates of \cite{[B-I]} are good 
for
our approach.

\hfill


\newpage
\begin{thebibliography}{References}
\bibitem{[B]}Bialy, M., {\em Convex Billiards and a theorem by
  E. Hopf}, Math. Z. {\bf 214} (1993) 147--154. 

\bibitem{[B-P1]}  Bialy, M., Polterovich, L., {\em Hopf type rigidity
  for Newton Equations},
  Math. Res. Lett. {\bf 2} (1995) 695--700.

\bibitem{[B-P2]} Bialy, M., Polterovich, L., {\em Hamiltonian
  diffeomorphisms and Lagrangian distributions}, GAFA 2(1992)173-210. 
 
\bibitem{[B-M]} Bialy, M., MacKay, R.S., {\em Variational properties
  of an elliptic nonlinear equation and rigidity}, preprint (1999).

\bibitem{[B-I]} Burago, D., Ivanov, S., {\em Riemannian tori without
  conjugate points are flat}, GAFA {\bf 4} (1994) 259--269. 

\bibitem{[C-F]} Croke, C., Fathi, A., {\em An inequality between energy
    and intersection}, Bull. Lond. Math. Soc. {\bf 22}, (1990) 489--494. 

\bibitem{[C-K]} Croke, C., Kleiner, B., {\em On tori without conjugate
    points}, Invent. Math. {\bf 120}, N2, (1995) 241--257. 

\bibitem{[C-I]} Contreras, G., Inturriaga, R., {\em Convex
    Hamiltonians without conjugate points}, preprint 1996.

\bibitem{[H]} Hopf, E., {\em Closed surfaces without conjugate
    points}, Proc. Nat. Acad. Sci. {\bf 34}, (1948) 47--51. 

\bibitem{[He]} Heber, J., {\em On the geodesic flow of tori without
    conjugate points}, Math. Z. {\bf 216}, N2, (1994) 209--216. 

\bibitem{[G]} Green, L., {\em A theorem of E. Hopf},
  Mich. Math. J. {\bf 5}, (1958) 31--34.

\bibitem{[Gi]} Ginzburg, V., {\em On closed trajectories of a charge
    in a magnetic field.  An application of symplectic geometry.}  In
    Thomas C.B. (editor) Contact and Symplectic geometry.
    Publ. Newton Institute {\bf 8}, (1996) 131--148.

\bibitem{[K]}  Knauf, A., {\em Closed orbits and converse KAM theory},
  Nonlinearity {\bf 3}, (1990) 961--973.

\end{thebibliography}
\end{document}